\documentclass{amsproc}
\usepackage{amssymb}
\usepackage[]{amsmath,amsthm,amssymb}
\usepackage[latin1]{inputenc}
\usepackage{graphicx}
\usepackage{amscd}

\theoremstyle{plain}

\newtheorem{definition}{Definition}

\newtheorem{remark}{Remark}

\newtheorem{theorem}{Theorem}
\numberwithin{equation}{section}

\input{tcilatex}

\begin{document}
\title[Orthogonality with Jacobi weight case ]{Hardy-type theorem for
orthogonal functions with respect to their zeros. The Jacobi weight case}
\author{L. D. Abreu$^{\star}$}
\address{$^{\star}$ Department of Mathematics, Faculty of Science and
Technology, University of Coimbra, P.Box 3008, 3001-454, Coimbra, Portugal}
\email{daniel@mat.uc.pt }
\urladdr{}
\thanks{}
\author{F. Marcellan$^{\dag}$}
\address{$^{\dag }$ Department of Mathematics, University Carlos III de
Madrid, Avenida de la Universidad 30, 28911 Legan{é}s, Spain}
\email{pacomarc@ing.uc3m.es}
\urladdr{}
\thanks{}
\author{S.B. Yakubovich$^{\ddag}$}
\address{$^{\ddag}$ Department of Pure Mathematics, Faculty of Sciences,
University of Porto, Campo Alegre st., 687, 4169-007 Porto, Portugal}
\email{syakubov@fc.up.pt}
\urladdr{}
\date{October 26, 2006}
\subjclass{44A15, 33C15, 33C10, 45E10, 42C15, 30D20}
\keywords{Zeros of special functions, Orthogonality, Jacobi weights, Mellin
transform on distributions, Entire functions, Bessel functions, Hyperbessel
functions}
\thanks{$^{\star }$ The work has been supported by CMUC and FCT
post-doctoral grant SFRH/BPD/26078/2005}
\thanks{$^{\dag }$ The work has been supported by Direcci{ó}n General de
Investigaci{ó}n, Ministerio de Educaci{ó}n y Ciencia of Spain, MTM 2006-
13000-C03-02 }
\thanks{$^{\ddag }$ The work has been supported, in part, by the "Centro de
Matem{á}tica" of the University of Porto}
\maketitle

\begin{abstract}
Motivated by G. H. Hardy's 1939 results \cite{Hardy} on functions orthogonal
with respect to their real zeros $\lambda _{n},\ n=1,2,\dots $, we will
consider, within the same general conditions imposed by Hardy, functions
satisfying an orthogonality with respect to their zeros with Jacobi weights
on the interval $(0,1)$, that is, the functions $f(z)=z^{\nu }F(z),\ \nu \in
\mathbb{R}$, where $F$ is entire and
\begin{equation*}
\int_{0}^{1}f(\lambda _{n}t)f(\lambda _{m}t)t^{\alpha }(1-t)^{\beta
}dt=0,\quad \alpha >-1-2\nu ,\ \beta >-1,
\end{equation*}%
when $n\neq m$. Considering all possible functions on this class we are lead
to the discovery of a new family of generalized Bessel functions including
Bessel and Hyperbessel functions as special cases.
\end{abstract}

\section{Introduction}

In his 1939 paper \cite{Hardy}, G. H. Hardy proved that, under certain
conditions, the only functions satisfying
\begin{equation*}
\int_{0}^{1}f(\lambda _{n}t)f(\lambda _{m}t)dt=0
\end{equation*}
if $m\neq n$, are the Bessel functions.

With a view to extend Hardy's results to the $q$-case, it was observed in
\cite{A} that a substantial part of Hardy's argument could be carried out
without virtually any change, when the measure $dt$ is replaced by an
arbitrary positive measure $du(t)$ on the interval $(0,1)$ and functions
satisfying
\begin{equation}
\int_{0}^{1}f(\lambda _{n}t)f(\lambda _{m}t)du(t)=0  \label{ortgen}
\end{equation}
if $m\neq n$, are considered.

In particular, the completeness of the set $\{f(\lambda _{n}t)\}$ holds for
such a general orthogonal system. Furthermore, there exists an associated
Lagrange-type sampling theorem for functions defined by integral transforms
whose kernel is defined by means of $f$ \cite{A}.

It is however not possible to identify, under such a degree of
generalization, the whole class of functions satisfying (\ref{ortgen}).
Therefore, we raise the question: Given a specific positive real measure on
the interval $(0,1)$, what are the resulting orthogonal functions? In the
case where this measure is the $d_{q}x$ arising from Jackson's $q$-integral,
it was shown in \cite{A} that the corresponding functions are the Jackson $q$%
-Bessel functions of the third type. In this note we will answer this
question in the case where $du(t)=t^{\alpha }(1-t)^{\beta }dt$ is the
measure associated with the orthogonal Jacobi polynomials on the interval $%
(0,1)$. That is, we will characterize in some sense the functions $f(z)=
z^\nu F(z)$, where $\nu \in \mathbb{R}$, $F$ is entire, which satisfy the
equality
\begin{equation}
\int_{0}^{1}f(\lambda _{n}t)f(\lambda _{m}t)t^{\alpha }(1-t)^{\beta }dt=0, \
\alpha > - 1 - 2\nu, \ \beta > -1,  \label{ortjacob}
\end{equation}
with membership of the classes $\mathcal{A}$ or $\mathcal{B}$ being
specified in the following definition (the same general restrictions are
imposed in \cite{Hardy} and \cite{A}).

\begin{definition}
Let $\nu, \alpha \in \mathbb{R}$ be such that $2\nu+\alpha > -1$ and $f(z)=
z^{\nu }F(z)$. The class $\mathcal{A}$ is constituted by entire functions $F$
of order less than two or of order two and minimal type having real positive
zeros $\lambda_n > 0, \ n=1,2, \ \dots$ with $\sum \lambda^{-2}_n < \infty$
situated symmetrically about the origin. The class $\mathcal{B}$, in turn,
contains entire functions $F(z)$ which have real but not necessarily
positive zeros $\lambda_n, \ n=1,2, \ \dots$ with $\sum \lambda^{-1}_n <
\infty$ and of order less than one or of order one and minimal type with $%
F(0)=1$.
\end{definition}

In this paper we will deal mainly with functions of the class $\mathcal{B}$.
Namely, we will show that these functions satisfy the Abel-type integral
equation of the second kind \cite{Yak} and will study their important
properties and series representations. However, in the last section we will
derive similar integral equation for the class $\mathcal{A}$ and will also
give its explicit solutions in terms of the series.

We will appeal in the sequel to the theory of the Mellin transform \cite{Tit}%
, \cite{Yak}, \cite{Zem}. As it is known the Mellin direct and inverse
transforms are defined by the formulas
\begin{equation*}
f^{*}(s)=\int_0^\infty f(x)x^{s-1}dx,\leqno(1.3)
\end{equation*}
\begin{equation*}
f(x)={\frac{1}{2\pi i}}\int_{\gamma-i\infty}^{\gamma+i\infty}
f^{*}(s)x^{-s}ds,\ s=\gamma+i\tau,\ x > 0,\leqno(1.4)
\end{equation*}
where integrals (1.3)- (1.4) exist as Lebesgue ones if we assume the
conditions $f \in L_1(\mathbb{R}_+; x^{\gamma-1}dx), \ f^* \in
L_1((\gamma-i\infty, \gamma +i\infty); d\tau)$, respectively. However
functions from the classes $\mathcal{A}$ or $\mathcal{B}$ can have a non-
integrable singularity at infinity and the classical Mellin transform (1.3)
generally does not exist. Therefore owing to \cite{Zem}, Ch. 4 we introduce
the Mellin transform for distributions from the dual space $%
M_{\gamma_1,\gamma_2}^\prime$ to Zemanian's testing -function space $%
M_{\gamma_1,\gamma_2}$ into the space of analytic functions in the open
vertical strip $\Omega_f:= \{ s \in \mathbb{C}: \gamma_1 < \mathrm{Re} s <
\gamma_2\}$ by the formula
\begin{equation*}
f^*(s):= \langle f(x), x^{s-1}\rangle, \ s \in \Omega_f.\leqno(1.5)
\end{equation*}
In this case the inversion integral (1.4) is convergent in $\mathcal{D}%
^\prime(\mathbb{R}_+)$, i.e. for any smooth function $\theta \in \mathcal{D}(%
\mathbb{R}_+) \subset M_{\gamma_1,\gamma_2}$ with compact support on $%
\mathbb{R}_+$ it holds
\begin{equation*}
\langle f(x), \theta(x)\rangle = \lim_{r\to \infty} \langle {\frac{1}{2\pi i}%
}\int_{\gamma-ir}^{\gamma+ir} f^{*}(s)x^{-s}ds, \ \theta (x) \rangle.\leqno%
(1.6)
\end{equation*}

\section{The Hardy type integral equation}

Suppose that $f\in \mathcal{B}$ and satisfies (\ref{ortjacob}). Let
\begin{equation*}
A_{n}=\int_{0}^{1}\left[ f(\lambda _{n}t)\right] ^{2}t^{\alpha }(1-t)^{\beta
}dt,\ \alpha >-1-2\nu ,\ \beta >-1.
\end{equation*}%
If $a_{n}(z)$ stands for the Fourier coefficients of the expansion of $f(zt)$
in terms of the orthonormal basis $\{A_{n}^{-\frac{1}{2}}f(\lambda _{n}t)\}$%
, then it is possible to show, following the arguments in \cite{Hardy}, that
\begin{equation*}
a_{n}(z)=\frac{1}{A_{n}^{1/2}}\int_{0}^{1}f(zt)f(\lambda _{n}t)t^{\alpha
}(1-t)^{\beta }dt=\frac{A_{n}f(z)}{f^{\prime }(\lambda _{n})(z-\lambda _{n})}%
\text{.}\leqno(2.1)
\end{equation*}%
In the sequel we will essentially continue to go along the lines of \cite%
{Hardy}, but since considerable number of adaptations have to be made in
order to deal with the presence of the Jacobi weight, we find it better to
provide the proof.

Using the Parseval identity for inner products we get
\begin{equation*}
\int_{0}^{1}f(zt)f(\zeta t)t^{\alpha }(1-t)^{\beta }dt=\sum_{n=0}^{\infty
}a_{n}(z)a_{n}(\zeta )\text{.}\leqno(2.2)
\end{equation*}

\begin{theorem}
Let $f \in \mathcal{B}$ and satisfy $(\ref{ortjacob})$ with $\alpha > -1 -
2\nu, \ \nu \in \mathbb{R}, \ \beta > -1$. Then it satisfies the integral
equation
\begin{equation*}
a\int_{0}^{z}u^{\nu +\alpha +1}(z-u)^{\beta }f(u)du=(az+1)\int_{0}^{z}u^{\nu
+\alpha }(z-u)^{\beta }f(u)du
\end{equation*}
\begin{equation*}
+z^{\nu +\alpha +\beta +1}f(z)A,\leqno(2.3)
\end{equation*}
where $a=F^{\prime }(0)$ and $A=- B(2\nu +\alpha +1,\beta +1)$, $B$ denotes
the Beta- function \cite{Erd}.
\end{theorem}

\begin{proof}
Substituting (2.1) into (2.2), we obtain, after some simplifications,
\begin{equation*}
\int_{0}^{1}f(zt)f(\zeta t)t^{\alpha }(1-t)^{\beta }dt=-f(z)f(\zeta )\frac{%
q(z)-q(\zeta )}{z-\zeta },\leqno(2.4)
\end{equation*}%
where
\begin{equation*}
q(z)=\sum_{n=1}^{\infty }\frac{A_{n}\lambda _{n}}{\{f^{\prime }(\lambda
_{n})\}^{2}}\left[ \frac{1}{z-\lambda _{n}}+\frac{1}{\lambda _{n}}\right] .
\end{equation*}%
Letting $\zeta \rightarrow 0$ in (2.4), the result is
\begin{equation*}
\int_{0}^{1}t^{\nu +\alpha }f(zt)(1-t)^{\beta }dt=-\frac{f(z)q(z)}{z}.
\end{equation*}%
The change of variables $u=zt$ gives
\begin{equation*}
\int_{0}^{z}u^{\nu +\alpha }(z-u)^{\beta }f(u)du=-z^{\nu +\alpha +\beta
}f(z)q(z).\leqno(2.5)
\end{equation*}%
Now, when $z$ is small enough, $f(z)\backsim z^{\nu }$ and $q(z)\sim
zq^{\prime }(0)$. Therefore, when $z\rightarrow 0$ we have
\begin{equation*}
\int_{0}^{z}u^{2\nu +\alpha }(z-u)^{\beta }du\sim -z^{2\nu +\alpha +\beta
+1}q^{\prime }(0)
\end{equation*}%
and, as a consequence,
\begin{equation*}
q^{\prime }(0)=-B(2\nu +\alpha +1,\beta +1).
\end{equation*}%
Rewriting (2.4) in the form
\begin{equation*}
\int_{0}^{1}t^{2\nu +\alpha }F(zt)F(\zeta t)(1-t)^{\beta }dt=-F(z)F(\zeta )%
\frac{q(z)-q(\zeta )}{z-\zeta }
\end{equation*}%
and differentiating in $\zeta $ we obtain
\begin{eqnarray*}
&&\int_{0}^{1}t^{2\nu +\alpha +1}F^{\prime }(\zeta t)F(zt)(1-t)^{\beta }dt \\
&=&-F(z)F^{\prime }(\zeta )\frac{q(z)-q(\zeta )}{z-\zeta }-F(z)F(\zeta )%
\frac{-q^{\prime }(\zeta )(z-\zeta )+q(z)-q(\zeta )}{(z-\zeta )^{2}}.
\end{eqnarray*}%
Setting $\zeta =0$, it gives
\begin{equation*}
F^{\prime }(0)\int_{0}^{1}t^{\nu +\alpha +1}(1-t)^{\beta }f(zt)dt=-{\frac{%
f(z)}{z^{2}}}\left[ q(z)(F^{\prime }(0)z+1)-Az\right] .
\end{equation*}%
So, if $a=F^{\prime }(0)$, the change of variable $zt=u$ yields
\begin{equation*}
a\int_{0}^{z}u^{\nu +\alpha +1}(z-u)^{\beta }f(u)du=-z^{\nu +\alpha +\beta
}f(z)\left[ (za+1)q(z)-Az\right] .
\end{equation*}%
Combining with (2.5) we get the integral equation (2.3).
\end{proof}

\section{The Hardy-type theorem}

In this section we will prove the Hardy type theorem \cite{Hardy}, which
will characterize all possible orthogonal functions with respect their zeros
from the class $\mathcal{B}$ under the Jacobi weight. Indeed, we have

\begin{theorem}
Let $f$ satisfy conditions of Theorem $1$. Then
\begin{equation*}
f(z)= \hbox{const.}\ z^\nu \sum_{n =0}^\infty a_{n}z^{n},\leqno(3.1)
\end{equation*}
where
\begin{equation*}
a_n = \frac{(a(\beta+1))^n}{\Gamma(\mu +n)}\prod_{j=1}^n \frac{(\mu)_j} {%
(\mu+ \beta + 1)_j - (\mu)_j},\ n=0, 1, \dots \ .\leqno(3.2)
\end{equation*}
Here $\mu= 2\nu+\alpha+1$, $\Gamma(z)$ is Euler's Gamma-function, $(b)_j$ is
the Pochhammer symbol \cite{Erd} and the empty product is equal to $1$.
Moreover, series $(3.1)$ represents an entire function of the order $\rho= {%
\frac{1}{\beta +2}} < 1$, when $\beta > -1$. In particular, the case $\beta
=0$ drives at the Hardy solutions in terms of the Bessel functions (cf. \cite%
{Hardy}, p.$43$).
\end{theorem}

\begin{proof}
Making again an elementary substitution $u=zt$ and cutting the multiplier $%
z^{\nu +\alpha +\beta +1}$ we write (2.3) in the form
\begin{equation*}
az\int_{0}^{1}t^{\nu +\alpha +1}(1-t)^{\beta
}f(zt)dt=(az+1)\int_{0}^{1}t^{\nu +\alpha }(1-t)^{\beta }f(zt)dt
\end{equation*}%
\begin{equation*}
-B(2\nu +\alpha +1,\beta +1)f(z).\leqno(3.3)
\end{equation*}%
Hence we observe that (3.3) is a second kind integral equation containing
two Erdélyi-Kober fractional integration operators with linear coefficients
(see \cite{Yak}, Ch. 3). Seeking a possible solution in terms of the series $%
z^{\nu }\sum_{n=0}^{\infty }a_{n}z^{n}$ where $a_{n}\neq 0,\ a_{0}=1$ and no
two consecutive $a_{n}$ vanish (see \cite{Hardy}, p. 43) we substitute it in
(3.3). After changing the order of integration and summation via the uniform
convergence, calculation the inner Beta-integrals, and elementary
substitutions we come out with the equality
\begin{equation*}
a\sum_{n=1}^{\infty }a_{n-1}z^{n+\nu }\left[ B(\mu +n,\beta +1)-B(\mu
+n-1,\beta +1)\right]
\end{equation*}%
\begin{equation*}
=\sum_{n=1}^{\infty }a_{n}z^{n+\nu }\left[ B(\mu +n,\beta +1)-B(\mu ,\beta
+1)\right] ,\ \mu =2\nu +\alpha +1.\leqno(3.4)
\end{equation*}%
Hence equating coefficients of the series and taking into account that $%
B(\mu +n,\beta +1)-B(\mu ,\beta +1)\neq 0,n\in \mathbb{N},\mu >0,\ \beta >-1$
we obtain the following recurrence relations
\begin{equation*}
a_{n}=a\ a_{n-1}\frac{B(\mu +n,\beta +1)-B(\mu +n-1,\beta +1)}{B(\mu
+n,\beta +1)-B(\mu ,\beta +1)}
\end{equation*}%
\begin{equation*}
=a^{n}\prod_{j=1}^{n}\frac{B(\mu +j,\beta +1)-B(\mu +j-1,\beta +1)}{B(\mu
+j,\beta +1)-B(\mu ,\beta +1)},\leqno(3.5)
\end{equation*}%
where $n\in \mathbb{N}_{0}$ and the empty product is equal to $1$. Moreover,
invoking the definition of the Pochhammer symbol $(b)_{k}={\frac{\Gamma (b+k)%
}{\Gamma (b)}}$, the reduction formula for Gamma-function $\Gamma
(z+1)=z\Gamma (z)$ and representing Beta-function via $B(b,c)={\frac{\Gamma
(b)\Gamma (c)}{\Gamma (b+c)}}$ we write (3.5) in the form (3.2).

Further, it is easily seen from (3.5) and asymptotic formula for the ratio
of Gamma-functions \cite{Erd} that
\begin{equation*}
\left|{\frac{a_n}{a_{n-1}}}\right|= \left|a\frac { B(\mu + n, \beta +1)-
B(\mu + n -1, \beta +1)}{B(\mu + n, \beta +1)- B(\mu,\beta +1)}\right|
\end{equation*}
\begin{equation*}
= \frac { \Gamma(\mu + n-1)}{\Gamma(\mu + n+\beta+ 1)} {\frac{|a|(\beta+1)}{%
B(\mu,\beta +1)- B(\mu + n, \beta +1)}}
\end{equation*}
\begin{equation*}
= O\left(n^{-(\beta+2)}\right) \to 0, \ n \to \infty, \ \beta > -1.
\end{equation*}
Thus $F(z)=\sum_{n =0}^\infty a_{n}z^{n}$ is an entire function. Let us
compute its order $\rho$. To do this we appeal to the familiar in calculus
Stolz's theorem. Hence we find
\begin{equation*}
\rho= \hbox{limsup}_{n \to \infty}{\frac{n \log n}{\log |a_n|^{-1}}} =
\lim_{n \to \infty}{\frac{(n+1) \log (n+1) - n \log n}{\log |a_{n+1}|^{-1}-
\log |a_{n}|^{-1}}}
\end{equation*}
\begin{equation*}
= \lim_{n \to \infty} {\frac{ \log n}{\log n^{\beta+2} }}= {\frac{1}{\beta +2%
}} < 1, \ \beta > -1.
\end{equation*}
Thus same arguments as in \cite{Hardy} guarantee that all our solutions
belong to the class $\mathcal{B}$. In particular, for $\beta =0$ we easily
get from (3.2) that
\begin{equation*}
a_n = \frac{a^n}{\Gamma(\mu +n)}\prod_{j=1}^n \frac{(\mu)_j} {(\mu + 1)_j -
(\mu)_j}= \frac{(a \mu)^n}{\Gamma(\mu +n)n!}, \ n=0, 1, \dots.
\end{equation*}
Then the series representation for the Bessel function \cite{Erd} yields
\begin{equation*}
f(z)= \hbox{const.}\ z^{-\alpha/2 }J_{\mu-1}\left(c z^{1/2}\right), \ c^2=
-4a\mu,
\end{equation*}
which slightly generalizes Hardy's solutions \cite{Hardy}, p.43. Theorem 2
is proved.
\end{proof}

\section{Properties of solutions and their particular cases}

First in this section we will apply the Mellin transform (1.5) to the
integral equation (2.3), reducing it to a certain functional equation. Then
we will solve this equation by methods of the calculus of finite differences
\cite{Mil} to obtain the value of the Mellin transform for solutions (3.1).

In fact, returning to (3.3) we consider $z=x\in \mathbb{R}_{+}$ and we apply
through the Mellin transform (1.5) taking into account the operational
formula
\begin{equation*}
\langle f(xt),x^{s-1}\rangle =t^{-s}f^{\ast }(s).
\end{equation*}%
Then using values of the elementary Beta-integrals we come out with the
following homogeneous functional equation
\begin{equation*}
aB(\nu +\alpha -s+1,\beta +1)f^{\ast }(s+1)=aB(\nu +\alpha -s,\beta
+1)f^{\ast }(s+1)
\end{equation*}%
\begin{equation*}
+B(\nu +\alpha -s+1,\beta +1)f^{\ast }(s)-B(2\nu +\alpha +1,\beta +1)f^{\ast
}(s),\leqno(4.1)
\end{equation*}%
where $s$ is a parameter of the Mellin transform (1.5) such that $s+1\in
\Omega _{f}$. Hence denoting by $h(s)=B(\nu +\alpha -s+1,\beta +1)$ and
invoking the condition $a\neq 0$ (see \cite{Hardy}, p.43) we rewrite (3.2)
as
\begin{equation*}
f^{\ast }(s+1)=H(s)f^{\ast }(s),\leqno(4.2)
\end{equation*}%
where the kernel $H(s)$ is given by
\begin{equation*}
H(s)={\frac{1}{a}}\frac{h(s)-h(-\nu )}{h(s)-h(s+1)}.\leqno(4.3)
\end{equation*}%
We can simplify $H(s)$ appealing to the basic properties for Beta-functions
\cite{Erd}. After straightforward calculations we get finally
\begin{equation*}
H(s)={\frac{1}{a(\beta +1)}}\frac{(s-\nu -\alpha )(h(s)-h(-\nu ))}{h(s)}.%
\leqno(4.4)
\end{equation*}%
Meanwhile,
\begin{equation*}
h(s)-h(-\nu )=\int_{0}^{1}(1-t)^{\beta }(t^{\nu +\alpha -s}-t^{2\nu +\alpha
})dt=\sum_{n=0}^{\infty }{\frac{(-\beta )_{n}}{n!}}
\end{equation*}%
\begin{equation*}
\times \left[ {\frac{1}{\nu +\alpha -s+n+1}}-{\frac{1}{2\nu +\alpha +n+1}}%
\right] =(\nu +s)\chi (s),
\end{equation*}%
where
\begin{equation*}
\chi (s)=\sum_{n=0}^{\infty }{\frac{(-\beta )_{n}}{n!(2\nu +\alpha +n+1)(\nu
+\alpha +n+1-s)}}.
\end{equation*}%
It is easily seen that $\chi (s)$ and $h(s)$ have the same simple poles at
the points $s=\nu +\alpha +n+1,\ n=0,1,\dots \ .$ As a consequence, taking
into account the expression of $h(s)$ in terms of Gamma-functions, we find
that $H(s)$ is meromorphic with simple poles $s=\nu +\alpha +\beta +n+2,\
n\in \mathbb{N}_{0}$. Furthermore, via the asymptotic behavior of the
Beta-function \cite{Erd}, we get $H(s)=O(s^{\beta +2}),\ s\rightarrow \infty
$.

The functional equation (4.2) can be formally solved by methods proposed in
\cite{Mil}, Ch. 11. So it has a general solution
\begin{equation*}
f^{\ast }(s)=\omega (s)\exp \left( \sum_{c}^{s}\log H(z)\Delta z\right) ,%
\leqno(4.5)
\end{equation*}%
where $\omega (s)$ is an arbitrary periodic function of $s$ of period $1$, $%
\sum_{c}^{s}$ means the operation of general summation \cite{Mil} with a
given constant $c$ and the function $\log H(s)$ is summable in this sense
(cf. \cite{Mil} Ch. 8). Clearly that it may be necessary to make suitable
cuts in the $s$ plane in view of the possible multi-valued nature of the
right-hand side of (4.5).

Let us consider interesting particular cases of the solutions (3.1) when $%
\beta =k\in \mathbb{N}_{0}$. Indeed, putting $\beta =k$ in (3.2) and using
the definition of Pochhammer's symbol we obtain
\begin{equation*}
{\frac{(\mu +k+1)_{j}}{(\mu )_{j}}}-1={\frac{\Gamma (\mu +k+j+1)}{\Gamma
(\mu +k+1)(\mu )_{j}}}-1={\frac{(\mu )_{k+j+1}}{(\mu )_{k+1}(\mu )_{j}}}-1
\end{equation*}%
\begin{equation*}
={\frac{(\mu +j)(\mu +j+1)\dots (\mu +k+j)}{(\mu )_{k+1}}}-1={\frac{jP_{k}(j)%
}{(\mu )_{k+1}}},
\end{equation*}%
where we denote by $P_{k}(j)=(j-\alpha _{1})\dots (j-\alpha _{k}),\ \alpha
_{i}\in \mathbb{C},\ i=1,\dots ,k$, a polynomial of degree $k$ with respect
to $j$ and evidently $P_{k}(j)>0,\ j\in \mathbb{N}$. Substituting the above
expression into (3.2) it becomes
\begin{equation*}
a_{n}=\frac{(a(k+1)(\mu )_{k+1})^{n}}{\Gamma (\mu +n)n!}\prod_{j=1}^{n}%
\prod_{i=1}^{k}\frac{1}{j-\alpha _{i}},\ n,k\in \mathbb{N}_{0}\leqno(4.6)
\end{equation*}%
and the empty products are equal to $1$. Consequently, solutions (3.1)
represent the so-called hyperbessel functions \cite{Yak}, Ch. 19 of the
order $k+1$, namely
\begin{equation*}
f(z)=\hbox{const.}\ z^{\nu }
\end{equation*}%
\begin{equation*}
\times {}_{0}F_{k+1}\left( 2\nu +\alpha +1,\ 1-\alpha _{1},\dots
,\ 1-\alpha _{k};\ a(k+1)(2\nu +\alpha +1)_{k+1}z\right) ,\ k\in
\mathbb{N}_{0}.\leqno(4.7)
\end{equation*}%
As it is known \cite{Yak} functions
$$F(z)={}_{0}F_{k+1}\left( 2\nu +\alpha +1,\ 1-\alpha _{1},\dots ,\
1-\alpha _{k};\ a(k+1)(2\nu +\alpha +1)_{k+1}z\right)$$
 in (4.7) satisfy the following $k+2$-th order linear differential equation
\begin{equation*}
{\frac{d}{dz}}\left( z{\frac{d}{dz}}+2\nu +\alpha \right)
\prod_{i=1}^{k}\left( z{\frac{d}{dz}}-\alpha _{i}\right)
F-a(k+1)(2\nu +\alpha +1 )_{k+1}F=0.\leqno(4.8)
\end{equation*}

\begin{remark}
{\rm We point out that in \cite{EM} the authors considered other
generalizations of the Bessel functions which also satisfy higher
order differential equations.}
\end{remark}

\begin{remark}
{\rm In particular, for $k=0$ we easily get Hardy's solutions in terms of the $%
{}_{0}F_{1}$ - functions and strictly prove their orthogonality
(1.2) by using the corresponding second order differential
equation (4.8). However, even for $k=1$ the direct proof of the
orthogonality property $(1.2)$ for solutions (4.7) is a difficult
task since we deal in this case with the third order differential
equation. As we aware, the orthogonality property for the
hyperbessel functions with respect to their zeros is yet unknown.
So, following conclusions of Theorem 2 we conjecture here this
fact as well as the orthogonality for all solutions (3.1) of the
class $\mathcal{B}$.}
\end{remark}

\section{Orthogonal functions of the class $\mathcal{A}$}

In the case $f\in \mathcal{A}$ we recall again the arguments in \cite{Hardy}
to write accordingly, the formula for Fourier coefficients $a_{n}(z)$ of $%
f(zt)$ (see (2.1)) as
\begin{equation*}
a_{n}(z)=\frac{2A_{n}\lambda _{n}}{f{\acute{}}(\lambda _{n})}\frac{f(z)}{%
z^{2}-\lambda _{n}^{2}}.\leqno(5.1)
\end{equation*}

We have

\begin{theorem}
If $f\in \mathcal{A}$ and satisfy $(\ref{ortjacob})$ with $\alpha >-1-2\nu
,\ \nu \in \mathbb{R},\ \beta >-1,$ then the integral equation holds
\begin{equation*}
a\int_{0}^{z}u^{\nu +\alpha +2}(z-u)^{\beta }f(u)du
\end{equation*}%
\begin{equation*}
=(az^{2}+2)\int_{0}^{z}u^{\nu +\alpha }(z-u)^{\beta }f(u)du+z^{\nu +\alpha
+\beta +1}f(z)A,\leqno(5.2)
\end{equation*}%
where $a=F^{\prime \prime }(0)$ and $A=-2B(2\nu +\alpha +1,\beta +1)$.
\end{theorem}

\begin{proof}
As in the proof of Theorem 1 we substitute (5.1) in (2.2) and we get the
equality
\begin{equation*}
\int_{0}^{1}f(zt)f(\zeta t)t^{\alpha }(1-t)^{\beta }dt=-f(z)f(\zeta )\frac{%
q(z)-q(\zeta )}{z^{2}-\zeta ^{2}},\leqno(5.3)
\end{equation*}%
where now
\begin{equation*}
q(z)=4\sum_{n=1}^{\infty }\frac{A_{n}\lambda _{n}^{2}}{\{f^{\prime }(\lambda
_{n})\}^{2}}\left[ \frac{1}{z^{2}-\lambda _{n}^{2}}+\frac{1}{\lambda _{n}^{2}%
}\right] .
\end{equation*}%
Letting $\zeta \rightarrow 0$ in (5.3), we obtain
\begin{equation*}
\int_{0}^{1}t^{\nu }f(zt)t^{\alpha }(1-t)^{\beta }dt=-\frac{f(z)q(z)}{z^{2}},
\end{equation*}%
which yields after the change $u=zt$
\begin{equation*}
\int_{0}^{z}u^{\nu +\alpha }(z-u)^{\beta }f(u)du=-z^{\nu +\alpha +\beta
-1}f(z)q(z).\leqno(5.4)
\end{equation*}%
When $z$ is small, $f(z)\backsim z^{\nu }$ and $q(z)\sim \frac{z^{2}}{2}%
q^{\prime \prime }(0)$. Therefore, as $z\rightarrow 0$ we get
\begin{equation*}
2\int_{0}^{z}u^{2\nu +\alpha }(z-u)^{\beta }du\sim -z^{2\nu +\alpha +\beta
+1}q^{\prime \prime }(0).
\end{equation*}%
Thus $q^{\prime \prime }(0)=-2B(2\nu +\alpha +1,\beta +1)=A$. If we rewrite
(5.3) in the form
\begin{equation*}
\int_{0}^{1}t^{2\nu +\alpha }F(zt)F(\zeta t)(1-t)^{\beta }dt=-F(z)F(\zeta )%
\frac{q(z)-q(\zeta )}{z^{2}-\zeta ^{2}}.
\end{equation*}%
then after differentiation with respect to $\zeta $ we get
\begin{eqnarray*}
&&\int_{0}^{1}t^{2\nu +\alpha +1}F^{\prime }(\zeta t)F(zt)(1-t)^{\beta }dt \\
&=&-F(z)F^{\prime }(\zeta )\frac{q(z)-q(\zeta )}{z^{2}-\zeta ^{2}}%
-F(z)F(\zeta )\frac{-q^{\prime }(\zeta )(z^{2}-\zeta ^{2})+2\zeta
(q(z)-q(\zeta ))}{(z^{2}-\zeta ^{2})^{2}}.
\end{eqnarray*}%
Differentiating again in $\zeta $ and setting $\zeta =0$, we deduce after
some simplifications
\begin{equation*}
F^{\prime \prime }(0)\int_{0}^{1}t^{\nu +\alpha +2}(1-t)^{\beta
}f(zt)dt=-f(z)\left[ \frac{F^{\prime \prime }(0)q(z)}{z^{2}}-\frac{A}{z^{2}}+%
\frac{2q(z)}{z^{4}}\right] .
\end{equation*}%
Letting $a=F^{\prime \prime }(0)$ and making the change $zt=u$ we obtain the
equation
\begin{equation*}
a\int_{0}^{z}u^{\nu +\alpha +2}(z-u)^{\beta }f(u)du=-z^{\nu +\alpha +\beta
-1}f(z)\left[ (z^{2}a+2)q(z)-Az^{2}\right] .
\end{equation*}%
Taking into account (5.4) we finally obtain the integral equation (5.2).
\end{proof}

An analog of the Hardy theorem for the class $\mathcal{A}$ is

\begin{theorem}
Let $f$ satisfy conditions of Theorem 3. Then
\begin{equation*}
f(z)= \hbox{const.} z^\nu\sum_{n=0}^\infty a_{2n}z^{2n}, \leqno(5.5)
\end{equation*}
where
\begin{equation*}
a_{2n}= \left({\frac{a(\beta+1)}{2}}\right)^n \prod_{j=1}^n \frac {(\beta
+2(\mu + 2j-1))(\mu)_{2(j-1)}}{(\mu+ \beta+1)_{2j}- (\mu)_{2j}}, \
n=0,1,\dots,
\end{equation*}
$\mu=2\nu+\alpha+1$ and the empty product is equal to $1$.
\end{theorem}

\begin{proof}
The substitution $u=zt$ and the value of $A$ reduce equation (5.2) to
\begin{equation*}
az^{2}\int_{0}^{1}t^{\nu +\alpha }(t^{2}-1)(1-t)^{\beta }f(zt)dt
\end{equation*}%
\begin{equation*}
=2\left[ \int_{0}^{1}t^{\nu +\alpha }(1-t)^{\beta }f(zt)dt-f(z)B(2\nu
+\alpha +1,\beta +1)\right] .\leqno(5.6)
\end{equation*}%
Since the case $f\in \mathcal{A}$ presumes the following series
representation
\begin{equation*}
f(z)=\hbox{const}.z^{\nu }\sum_{n=0}^{\infty }a_{2n}z^{2n},
\end{equation*}%
with $a_{2n}\neq 0$ for any $n\in \mathbb{N}_{0}$, we substitute it into
(5.6), change the order of integration and summation, and calculate the
inner Beta-integrals. Therefore denoting by $\mu =2\nu +\alpha +1$ we obtain
\begin{equation*}
a\sum_{n=1}^{\infty }a_{2(n-1)}z^{2n}\left[ B(\mu +2n,\beta +1)-B(\mu
+2(n-1),\beta +1)\right]
\end{equation*}%
\begin{equation*}
=2\sum_{n=1}^{\infty }a_{2n}z^{2n}\left[ B(\mu +2n,\beta +1)\ -\ B(\mu
,\beta +1)\right] .
\end{equation*}%
Hence equating coefficients of the series the following recurrence relations
appear
\begin{equation*}
a_{2n}={\frac{a}{2}}a_{2(n-1)}\frac{B(\mu +2n,\beta +1)-B(\mu +2(n-1),\beta
+1)}{B(\mu +2n,\beta +1)\ -\ B(\mu ,\beta +1)}
\end{equation*}%
\begin{equation*}
=\left( {\frac{a}{2}}\right) ^{n}\prod_{j=1}^{n}\frac{B(\mu +2j,\beta
+1)-B(\mu +2(j-1),\beta +1)}{B(\mu +2j,\beta +1)-B(\mu ,\beta +1)},\leqno%
(5.7)
\end{equation*}%
where $\ n=0,1,2,\dots $, and the empty product is equal to $1$. In the same
way as (3.2) it can be written in the form
\begin{equation*}
a_{2n}=\left( {\frac{a(\beta +1)}{2}}\right) ^{n}\prod_{j=1}^{n}\frac{(\beta
+2(\mu +2j-1))(\mu )_{2(j-1)}}{(\mu +\beta +1)_{2j}-(\mu ){2j}},\
n=0,1,\dots \ .
\end{equation*}
\end{proof}

Putting $\beta =0$ in (5.7) by straightforward calculations functions (5.5)
become
\begin{equation*}
f(z)=\hbox{const.}z^{{\frac{1-\alpha }{2}}}J_{{\frac{\mu }{2}}-1}\left(
(-a\mu )^{1/2}z\right) ,
\end{equation*}%
which slightly generalizes Hardy's solutions \cite{Hardy}, p.42.

It is easily seen that the series in (5.5) is an entire function. Finally we
compute the order of these solutions. By similar arguments as in Theorem 2
we have
\begin{equation*}
\rho =\hbox{limsup}_{n\rightarrow \infty }{\frac{2n\log 2n}{\log
|a_{2n}|^{-1}}}=\lim_{n\rightarrow \infty }{\frac{2(n+1)\log 2(n+1)-2n\log 2n%
}{\log |a_{2(n+1)}|^{-1}-\log |a_{2n}|^{-1}}}
\end{equation*}%
\begin{equation*}
=\lim_{n\rightarrow \infty }{\frac{2\log n}{\log n^{\beta +2}}}={\frac{2}{%
\beta +2}}<2,\ \beta >-1.
\end{equation*}%
Therefore all solutions belong to the class $\mathcal{A}$.

\begin{remark}
{\rm Our final conjecture is that solutions $f\in \mathcal{A}$ of
the Hardy-type integral equation (5.2) are orthogonal with respect
to their own zeros for the Jacobi weight}.
\end{remark}

\end{document}